%&amstex
\input amsppt.sty
\hsize=12.5cm
\vsize=540pt
\baselineskip=16truept
\NoRunningHeads
\TagsOnRight
%-----------------------------------------
 \define\CC{\Bbb C}
 \define\CF{\Cal F}
 \define\CH{\Cal H}
 \define\CO{\Cal O}
 \redefine\Delta{\varDelta}
 \define\eps{\varepsilon}
 \define\halfskip{\vskip6pt}
 \define\intt{\operatorname{int}}
 \redefine\kappa{\varkappa}
 \font\msbm=msbm10
 \define\MX{\text{\msbm X}}
 \define\NN{\Bbb N}
 \define\pder#1#2{\frac{\partial#1}{\partial#2}}
 \redefine\phi{\varphi}
 \define\proj{\operatorname{pr}}
 \define\PSH{\Cal P\Cal S\Cal H}
 
 \define\Reg{\operatorname{Reg}}
 \define\Sing{\operatorname{Sing}}
 
 \define\tdots{\times\dots\times}
 \define\too{\longrightarrow}
 \define\wdht{\widehat}
 \define\wdtl{\widetilde}
%--------------------
\catcode`@=11
\def\qed{\ifhmode\textqed\else\ifinner\quad\square
   \else\eqno\square\fi\fi}
\def\textqed{{\unskip\nobreak\penalty50
    \quad\hbox{}\nobreak\hfil$\m@th\square$
    \parfillskip=0pt \finalhyphendemerits=0\par}}
\catcode`@=\active
%==================================

\document
\topmatter
\title
An extension theorem for separately holomorphic functions with singularities
\endtitle
\author
Marek Jarnicki (Krak\'ow)
\footnote"${}^{(\dag)}$"
{\;\hbox{Research partially supported by the KBN grant No\. 5 P03A 033 21.}}
Peter Pflug (Oldenburg)
\footnote"${}^{(\ddag)}$"
{\;Research partially supported by the Nieders\"achsisches
Ministerium f\"ur Wissenschaft und
\hbox{\hskip12pt Kultur, Az. 15.3 -- 50 113(55) PL.}}
\endauthor
\address
Uniwersytet Jagiello\'nski \newline
Instytut Matematyki \newline
30-059 Krak\'ow, Reymonta 4, Poland\newline
{\it E-mail: }{\rm jarnicki{\@}im.uj.edu.pl}
\endaddress

\address
Carl von Ossietzky Universit\"at Oldenburg
\newline Fachbereich Mathematik\newline
Postfach 2503\newline
D-26111 Oldenburg, Germany\newline
{\it E-mail: }{\rm pflug{\@}mathematik.uni{-oldenburg.de}}
\endaddress

\abstract
Let $D_j\subset\CC^{k_j}$ be a pseudoconvex domain and let $A_j\subset D_j$
be a locally pluripolar set, $j=1,\dots,N$. Put
$$
X:=\bigcup_{j=1}^N A_1\tdots A_{j-1}\times D_j\times A_{j+1}\tdots A_N
\subset\CC^{k_1+\dots+k_N}.
$$
Let $U$ be an open connected neighborhood of $X$ and let $M\varsubsetneq U$
be an analytic subset. Then there
exists an analytic subset $\wdht M$ of the `envelope of holomorphy' $\wdht X$
of $X$ with $\wdht M\cap X\subset M$ such that for every function $f$
separately holomorphic on $X\setminus M$ there exists an $\wdht f$ holomorphic
on $\wdht X\setminus\wdht M$ with $\wdht f|_{X\setminus M}=f$.
The result generalizes special cases which were
studied in \cite{\"Okt 1998}, \cite{\"Okt  1999}, \cite{Sic 2000}, and
\cite{Jar-Pfl 2001}.
\endabstract


\endtopmatter

\noindent{\bf 1. Introduction. Main Theorem.} Let $N\in\NN$, $N\geq2$, and
let
$$
\varnothing\neq A_j\subset D_j\subset\CC^{k_j},
$$
where $D_j$ is a domain, $j=1,\dots,N$. We define an {\it $N$--fold cross}\
$$
\align
X&:=\MX(A_1,\dots,A_N;D_1,\dots,D_N)\\
&:=\bigcup_{j=1}^N
A_1\tdots A_{j-1}\times D_j\times A_{j+1}\tdots A_N
\subset\CC^{k_1+\dots+k_N}.\tag1
\endalign
$$
Observe that $X$ is connected.

\halfskip

Let $\Omega\subset\CC^n$ be an open set and let $A\subset\Omega$. Put
$$
h_{A,\Omega}:=\sup\{u\:\; u\in\PSH(\Omega),\;u\leq1 \text{ on }
\Omega,\; u\leq0 \text{ on } A\},
$$
where $\PSH(\Omega)$ denotes the set of all functions plurisubharmonic on
$\Omega$. Define
$$
\omega_{A,\Omega}:=\lim_{k\to+\infty}h^\ast_{A\cap\Omega_k,\Omega_k},
$$
where $(\Omega_k)_{k=1}^\infty$ is a sequence of relatively compact open
sets $\Omega_k\subset\Omega_{k+1}\subset\subset\Omega$ with
$\bigcup_{k=1}^\infty\Omega_k=\Omega$ ($h^\ast$ denotes the upper
semicontinuous regularization of $h$). Observe that the definition is
independent of the chosen exhausting sequence $(\Omega_k)_{k=1}^\infty$.
Moreover, $\omega_{A,\Omega}\in\PSH(\Omega)$.

\halfskip

For an $N$--fold cross $X=\MX(A_1,\dots,A_N;D_1,\dots,D_N)$ put
$$
\wdht X:=\{(z_1,\dots,z_N)\in D_1\tdots D_N\:
\sum_{j=1}^N\omega_{A_j,D_j}(z_j)<1\};
$$
notice that $\wdht X$ may be empty.
Observe that $\wdht X$ is pseudoconvex if $D_1,\dots,D_N$ are pseudoconvex
domains.

\halfskip

We say that a subset $\varnothing\neq A\subset\CC^n$ is {\it locally pluriregular}
if $h^\ast_{A\cap\Omega,\Omega}(a)=0$ for any $a\in A$ and for any open
neighborhood $\Omega$ of $a$ (in particular, $A\cap\Omega$ is non-pluripolar).

Note that if $A_1,\dots,A_N$ are locally pluriregular, then
$X\subset\wdht X$ and, moreover, $\wdht X$ is connected (Lemma 5).

\halfskip

Let $U$ be a connected neighborhood of $X$ and let
$M\varsubsetneq U$ be an analytic subset ($M$ may be empty).
We say that a function $f\:X\setminus M\too\CC$ is {\it separately holomorphic}
($f\in\CO_s(X\setminus M)$) if for any $(a_1,\dots,a_N)\in A_1\tdots A_N$
and $j\in\{1,\dots,N\}$ the function
$f(a_1,\dots,a_{j-1},\cdot,a_{j+1},\dots,a_N)$ is holomorphic in the domain
$\{z_j\in D_j\:(a_1,\dots,a_{j-1},z_j,a_{j+1},\dots,a_N)\not\in M\}$.

\halfskip

The main result of our paper is the following extension theorem for
separately holomorphic functions. \footnote"${}^{(\ddag)}$"{\; We like to
thank Professor J\'ozef Siciak for turning our attention to the problem.}

\proclaim{Main Theorem} Let $D_j\subset\CC^{k_j}$ be a pseudoconvex domain
and let $A_j\subset D_j$ be a locally pluriregular set, $j=1,\dots,N$.
Let $M\varsubsetneq U$ be an analytic subset of an open connected neighborhood
$U$ of $X=\MX(A_1,\dots,A_N;D_1,\dots,D_N)$ ($M$ may be empty).
Then there exists a pure one--codimensional analytic subset
$\wdht M\subset\wdht X$ such that:

$\bullet$ $\wdht M\cap U_0\subset M$ for an open neighborhood $U_0$ of $X$,
$U_0\subset U$,

$\bullet$ for every $f\in\CO_s(X\setminus M)$ there exists exactly one
$\wdht f\in\CO(\wdht X\setminus\wdht M)$ with $\wdht f|_{X\setminus M}=f$.

Moreover, if $U=\wdht X$, then we can take $\wdht M:=$ the union of all
one--codimensional irreducible components of $M$.
\endproclaim

The proof will be given in Sections 3 (the case $U=\wdht X$) and 4
(the general case).

\proclaim{Remark} Notice that the Main Theorem may be generalized to the case
where $D_j$ is Riemann--Stein domains over $\CC^{k_j}$, $j=1,\dots,N$.
\endproclaim

Observe that in the case $M=\varnothing$, $N=2$, the Main Theorem is nothing
else than the following cross theorem.

\proclaim\nofrills{Theorem 1}\;{\rm (cf\. \cite{Ale-Zer 2001}).}
Let $D\subset\CC^p$, $G\subset\CC^q$ be pseudoconvex
domains and let $A\subset D$, $B\subset G$ be locally pluriregular.
Put $X:=\MX(A,B;D,G)$.
Then for any $f\in\CO_s(X)$ there exists exactly one
$\wdht f\in\CO(\wdht X)$ with $\wdht f=f$ on $X$.
\endproclaim

\proclaim{Remark}\rm (a) $M=\varnothing$: There is a long list of papers
discussing  the case $N=2$ (under various assumptions): \cite{Sic 1969},
\cite{Zah 1976}, \cite{Sic 1981}, \cite{Shi 1989}, \cite{Ngu-Zer 1991},
\cite{Ngu 1997}, \cite{Ale-Zer 2001}.
The case $N\geq2$, $k_1=\dots=k_N=1$ can be found in \cite{Sic 1981}.
The general case $N\geq2$, $k_1,\dots,k_N\geq1$ was solved in
\cite{Ngu-Zer 1995}\;
\footnote"${}^{(\ddag)}$"{\;We like to thank Professor Nguyen Thanh Van for
calling our attention to that paper.}.

\halfskip

\noindent(b) $M\neq\varnothing$: J\. Siciak \cite{Sic 2000} solved the
case: $N\geq2$, $k_1=\dots=k_N=1$,
$D_1=\dots=D_N=\CC$, $M=P^{-1}(0)$, where $P$ is a non-zero polynomial of
$N$ complex variables. The special subcase $N=2$, $P(z,w):=z-w$ had been
studied in \cite{\"Okt 1998}. The general case for $N=2$, $k_1=k_2=1$
was solved in \cite{Jar-Pfl 2001}; see also \cite{\"Okt 1999} for a partial
discussion of the case $N=2$, $k_1,k_2\geq1$.
\endproclaim

It remains an open problem what happens if the singular set $M$ is assumed
to be, for instance, a pluripolar relatively closed subset of $U$ (see also
\cite{Chi-Sad 1988}). The following particular case will be proved in
Section 5.

\proclaim{Theorem 2} Let $A, B\subset\CC$ be locally regular.
Let $M\subset\CC^2$ be a closed pluripolar set and let
$X:=\MX(A,B;\CC,\CC)$.
Then there exist a closed pluripolar set $\wdht M\subset\CC^2$ and
$A'\subset A$, $B'\subset B$ with $A\setminus A'$, $B\setminus B'$ being
polar, such that for every $f\in\CO_s(X\setminus M)$ there exists exactly one
$\wdht f\in\CO(\CC^2\setminus\wdht M)$ with
$\wdht f|_{X'\setminus(M\cup\wdht M)}=f|_{X'\setminus(M\cup\wdht M)}$,
where $X':=\MX(A',B';\CC,\CC)$.
\endproclaim

Notice that Theorem 2 may be thought as a partial generalization of
\cite{Sic 2000}.

\halfskip
\noindent{\bf 2. Auxiliary results.}
The following lemma gathers a few standard results, which will be
frequently used in the sequel.

\proclaim\nofrills{Lemma 3}\;
{\rm (cf\. \cite{Kli 1991}, \cite{Jar-Pfl 2000}, \S\;3.5).}
{\rm (a)} Let $\Omega\subset\CC^n$ be a bounded open set and let
$A\subset\Omega$. Then:

$\bullet$ If $P\subset\CC^n$ is pluripolar, then
$h^\ast_{A\setminus P,\Omega}=h^\ast_{A,\Omega}$.

$\bullet$ $h_{A_k\cap\Omega_k,\Omega_k}^\ast\searrow h_{A,\Omega}^\ast$
(pointwise on $\Omega$) for any sequence of open sets
$\Omega_k\nearrow\Omega$ and any sequence $A_k\nearrow A$.

$\bullet$ $\omega_{A,\Omega}=h^\ast_{A,\Omega}$.

$\bullet$ The following conditions are equivalent:

\hang{for any connected component $S$ of $\Omega$ the set $A\cap S$ is
non-pluripolar;}

\hang{$h^\ast_{A,\Omega}(z)<1$ for any $z\in\Omega$.}

$\bullet$ If $A$ is non-pluripolar, $0<\alpha<1$, and
$\Omega_\alpha:=\{z\in\Omega\: h^\ast_{A,\Omega}(z)<\alpha\}$,
then for any connected component $S$ of $\Omega_\alpha$ the set
$A\cap S$ is non-pluripolar (in particular, $A\cap S\neq\varnothing$).

$\bullet$ If $A$ is locally pluriregular, $0<\alpha\leq1$, and
$\Omega_\alpha:=\{z\in\Omega\: h^\ast_{A,\Omega}(z)<\alpha\}$, then
$h^\ast_{A,\Omega}=\alpha h^\ast_{A,\Omega_\alpha}$ on $\Omega_\alpha$.

\halfskip

\noindent{\rm (b)} Let $\Omega\subset\CC^n$ be an open set and let
$A\subset\Omega$. Then:

$\bullet$ $\omega_{A,\Omega}\in\PSH(\Omega)$.

$\bullet$ If $A$ is locally pluriregular, then $\omega_{A,\Omega}(a)=0$
for any $a\in A$.

$\bullet$ If $P\subset\CC^n$ is pluripolar, then
$\omega_{A\setminus P,\Omega}=\omega_{A,\Omega}$.

$\bullet$ If $A$ is locally pluriregular and $P\subset\CC^n$ is pluripolar,
then $A\setminus P$ is locally pluriregular.
\endproclaim

Moreover, we get:

\proclaim{Lemma 4} {\rm (a)} Let $A_j\subset\CC^{k_j}$ be locally pluriregular,
$j=1,\dots,N$, then $A_1\tdots A_N$ is locally pluriregular.

\noindent{\rm (b)} Let $A_j\subset\Omega_j\Subset\CC^{k_j}$,
$\Omega_j$ a domain, $A_j$ locally pluriregular, $j=1,\dots,N$, $N\geq2$.
Put
$$
\Omega:=\{(z_1,\dots,z_N)\in\Omega_1\tdots\Omega_N\:
\sum_{j=1}^Nh_{A_j,\Omega_j}^\ast(z_j)<1\}
$$
(observe that $A_1\tdots A_N\subset\Omega$). Then
$$
h_{A_1\tdots A_N, \Omega}^\ast=\sum_{j=1}^N h_{A_j,\Omega_j}^\ast
\;\text{ on } \Omega.
$$
\endproclaim

\demo{Proof} (a) is an immediate consequence of the product property
for the relatively extremal function
$$
h_{B_1\tdots B_N, \Omega_1\tdots\Omega_N}^\ast=
\max\{h_{B_j,\Omega_j}^\ast\: j=1,\dots,N\};
$$
cf\. \cite{Ngu-Sic 1991}.

\halfskip

(b) First observe that
$$
\sum_{j=1}^N h_{A_j,\Omega_j}^\ast\leq h_{A_1\tdots A_N, \Omega}^\ast
\;\text{ on } \Omega.
$$
To get the opposite inequality we proceed by induction on
$N\geq2$.

Let $N=2$: The proof of this step is taken from \cite{Sic 1981}.
For the reader's convenience we repeat the details.

Put $u:=h_{A_1\times A_2,\Omega}^\ast\in\PSH(\Omega)$ and fix a point
$(a_1,a_2)\in\Omega$.

If $a_1\in A_1$ (thus $h_{A_1,\Omega_1}^\ast(a_1)=0$), then
$u(a_1,\cdot)\in\PSH(\Omega_2)$ with $u(a_1,\cdot)\leq1$ and
$u(a_1,\cdot)\leq0$ on $A_2$. Therefore,
$$
u(a_1,\cdot)\leq h_{A_2,\Omega_2}^\ast=h_{A_1,\Omega_1}^\ast(a_1)+
h_{A_2,\Omega_2}^\ast\;\text{ on } \Omega_2.
$$
In particular, $u(a_1,a_2)\leq h_{A_1,\Omega_1}^\ast(a_1)+
h_{A_2,\Omega_2}^\ast(a_2)$.

Observe that the same argument shows that if $a_2\in A_2$, then
$u(\cdot,a_2)\leq h_{A_1,\Omega_1}^\ast$ on $\Omega_1$.

If $a_1\notin A_1$, then $h_{A_1,\Omega_1}^\ast(a_1)+
h_{A_2,\Omega_2}^\ast(a_2)<1$ and therefore
$
\alpha:=1-h_{A_2,\Omega_2}^\ast(a_2)\in(0,1]$.
Put
$$
(\Omega_2)_\alpha:=\{z_2\in\Omega_2\: h_{A_2,\Omega_2}^\ast(z_2)<\alpha\}.
$$
It is clear that $A_2\subset(\Omega_2)_\alpha\ni a_2$. Put
$$
v:=\frac1\alpha\Big(u(a_1,\cdot)-h_{A_1,\Omega_1}^\ast(a_1)\Big)
\in\PSH((\Omega_2)_\alpha).
$$
Then $v\leq1$ and $v\leq0$ on $A_2$. Therefore, by Lemma 3(a),
$$
v\leq h_{A_2,(\Omega_2)_\alpha}^\ast(a_2)=\frac1\alpha
h_{A_2,\Omega_2}^\ast(a_2)\; \text{ on } (\Omega_2)_\alpha.
$$
Consequently, $u(a_1,a_2)\leq h_{A_1,\Omega_1}^\ast(a_1)+
h_{A_2,\Omega_2}^\ast(a_2)$, which finishes the proof for $N=2$.

\halfskip

Now, assume that the formula is true for $N-1\geq2$. Put
$$
\wdtl\Omega:=\{(z_1,\dots,z_{N-1})\in\Omega_1\tdots\Omega_{N-1}\:
\sum_{j=1}^{N-1}h_{A_j,\Omega_j}^\ast(z_j)<1\}
$$
and fix an arbitrary $z=(\wdtl z,z_N)\in\Omega$. Obviously,
$\wdtl z\in\wdtl\Omega$. In virtue of the inductive hypothesis, we conclude
that
$$
h_{A_1\tdots A_{N-1},\wdtl\Omega}^\ast(\wdtl z)=
\sum_{j=1}^{N-1} h_{A_j,\Omega_j}^\ast(z_j).\tag2
$$
Now we apply the case $N=2$ to the following situation:
$$
\Omega':=\{(\wdtl w,w_N)\in\wdtl\Omega\times\Omega_N\:
h_{A_1\tdots A_{N-1},\wdtl\Omega}^\ast(\wdtl w)+h_{A_N,\Omega_N}
^\ast(w_N)<1\}.
$$
So
$$
h_{A_1\tdots A_{N-1},\wdtl\Omega}^\ast(\wdtl w)+h_{A_N,\Omega_N}
^\ast(w_N)=h_{A_1\tdots A_N, \Omega'}^\ast(\wdtl w,w_N),\quad
(\wdtl w,w_N)\in\Omega'.
$$
Note that $\Omega'=\Omega$. Hence
$$
h_{A_1\tdots A_N, \Omega}^\ast(\wdtl z,z_N)=
h_{A_1\tdots A_{N-1},\wdtl\Omega}^\ast(\wdtl z)+h_{A_N,\Omega_N}
(z_N)
\overset{(2)}\to=\sum_{j=1}^N h_{A_j,\Omega_j}^\ast(z_j).
$$
\qed\enddemo

\proclaim{Lemma 5} Let $X=\MX(A_1,\dots,A_N;D_1,\dots, D_N)$ be an
$N$--fold cross as in {\rm (1)}. If $A_1,\dots, A_N$ are locally pluriregular,
then $\wdht X$ is a domain.
\endproclaim

\demo{Proof} Using exhaustion by bounded domains we may assume that
the $D_j$'s are bounded.

We know that $X\subset\wdht X$. Let $z^0=(z_1^0,\dots,z_N^0)\in\wdht X$
be an arbitrary point.

If $\sum_{j=2}^N h_{A_j,\Omega_j}^\ast(z_j^0)=0$, then $D_1\times
\{(z_2^0,\dots,z_N^0)\}\subset\wdht X$. Therefore, $z^0$ can be joined
inside $D_1\times\{(z_2^0,\dots,z_N^0)\}$
with $(a_1,z_2^0,\dots,z_N^0)$ for some $a_1\in A_1$.

If $\sum_{j=2}^N h_{A_j,\Omega_j}^\ast(z_j^0)=:\eps>0$, put
$$
(D_1)_{1-\eps}:=\{z_1\in D_1\: h_{A_1,D_1}^\ast(z_1)<1-\eps\}.
$$
Then, in virtue of Lemma 3(a), the connected component $S$ of
$(D_1)_{1-\eps}$, that contains $z_1^0$, intersects $A_1$. Therefore,
$z^0$ can be joined inside $S\times\{(z_2^0,\dots,z_N^0)\}\subset\wdht X$
with $(a_1,z_2^0,\dots,z_N^0)$ for some $a_1\in A_1$.

Now we repeat the above argument for the second component of the point
$(a_1, z_2^0, \dots, z_N^0)$. Finally, the point $z^0$ can be joined
inside $\wdht X$ with $(a_1,\dots,a_N)\in A_1\tdots A_N
\subset X$. Since $X$ is connected, the proof is completed.
\qed\enddemo

\proclaim\nofrills{Lemma 6}\; {\rm (Identity theorems).}
{\rm (a)} Let $\Omega\subset\CC^n$ be a domain and let $A\subset\Omega$
be non-pluripolar. Then any $f\in\CO(\Omega)$ with $f|_A=0$ vanishes
identically on $\Omega$.

\noindent{\rm (b)}
Let $D\subset\CC^p, G\subset\CC^q$ be domains, let $A\subset D$,
$B\subset G$ be locally pluriregular sets, and let $X:=\MX(A,B;D,G)$.
Let $M\varsubsetneq U$ be an analytic subset of an open connected
neighborhood $U$ of $X$. Assume that $A'\subset A$, $B'\subset B$
are such that:

$\bullet$ $A\setminus A'$ and
$B\setminus B'$ are pluripolar (in particular, $A'$, $B'$ are also locally
pluriregular),

$\bullet$ $M_z:=\{w\in G\: (z,w)\in M\}\neq G$ for any $z\in A'$,

$\bullet$ $M^w:=\{z\in D\: (z,w)\in M\}\neq D$ for any $w\in B'$.

\noindent Then:

\noindent {\rm (b${}_1$)} If $f\in\CO_s(X\setminus M)$ and
$f=0$ on $A'\times B'\setminus M$\; \footnote"${}^{(\ddag)}$"{\; Here and in the sequel
to simplify notation we write $P\times Q\setminus M$ instead of $(P\times Q)
\setminus M$.}, then $f=0$ on $X\setminus M$.

\noindent {\rm (b${}_2$)} If $g\in\CO(U\setminus M)$ and $g=0$ on
$A'\times B'\setminus M$, then $g=0$ on $U\setminus M$.
\endproclaim

\demo{Proof} (a) is obvious.

\halfskip

(b${}_1$)
Take a point $(a_0,b_0)\in X\setminus M$. We may assume that
$a_0\in A$. Since $A\setminus A'$ is pluripolar, there exists a sequence
$(a_k)_{k=1}^\infty\subset A'$ such that $a_k\too a_0$. The set
$Q:=\bigcup_{k=0}^\infty M_{a_k}$ is pluripolar. Consequently, the set
$B'':=B'\setminus Q$ is non-pluripolar. We have $f(a_k,w)=0$, $w\in B''$,
$k=1,2,\dots$. Hence $f(a_0,w)=0$ for any $w\in B''$. Finally,
$f(a_0,w)=0$ on $G\setminus M_{a_0}\ni b_0$.

\halfskip

(b${}_2$) Take an $a_0\in A'$. Since $M_{a_0}\neq G$, there exists a
$b_0\in B'\setminus M_{a_0}$. Let
$P=\Delta_{a_0}(r)\times\Delta_{b_0}(r)\subset U\setminus M$
($\Delta_{z_0}(r)\subset\CC^p$ denotes the polydisc
with center $z_0\in\CC^p$ and radius $r>0$). Then
$g(\cdot,w)=0$ on $A'\cap\Delta_{a_0}(r)$ for any $w\in B'\cap\Delta_{b_0}(r)$.
The set $A'\cap\Delta_{a_0}(r)$ is non-pluripolar.
Hence $g(\cdot,w)=0$ on $\Delta_{a_0}(r)$ for any $w\in B'\cap\Delta_{b_0}(r)$.
By the same argument for the second variable we get $g=0$ on $P$ and,
consequently, on $U$.
\qed\enddemo

\halfskip

\noindent{\bf 3. Proof of the Main Theorem in the case where $\bold{U=\wdht X}$.}
We proceed by several reduction
steps. First observe that, by Lemma 6(a), the function $\wdht f$ is uniquely
determined (if exists).

\proclaim{Step 1} To prove the Main Theorem for $M\neq\varnothing$
it suffices to consider only the case where $M$ is pure one--codimensional.
\endproclaim

\demo{Proof of Step 1} Since $\wdht X$ is pseudoconvex, the arbitrary
analytic set $M\subset\wdht X$ can be written as
$$
M=\{z\in\wdht X\: g_1(z)=\dots=g_k(z)=0\},
$$
where $g_j\in\CO(\wdht X)$, $g_j\not\equiv0$, $j=1,\dots,k$. Then
$M_j:=g_j^{-1}(0)$ is pure one--codimensional.

Take an $f\in\CO_s(X\setminus M)$.
Observe that $f_j:=f|_{X\setminus M_j}\in\CO_s(X\setminus M_j)$. By the
reduction assumption there exists an $\wdht f_j\in\CO(\wdht X\setminus M_j)$
such that $\wdht f_j=f$ on $X\setminus M_j$.
In virtue of Lemma 6(a), gluing the functions $(\wdht f_j)_{j=1}^k$,
leads to an $\wdht f\in\CO(\wdht X\setminus M)$ with $\wdht f=\wdht f_j$
on $\wdht X\setminus M_j$, $j=1,\dots,k$. Therefore, $\wdht f=f$ on
$X\setminus M$.

Finally, since codim$(M\setminus\wdht M)\geq2$, the function $\wdht f$
extends holomorphically across $M\setminus\wdht M$.
\qed\enddemo

{\it From now on we assume that $M$ is empty or pure one--codimensional.}

\proclaim{Step 2} To prove the Main Theorem it suffices to consider only the
case where $M$ is empty or pure one--codimensional and $D_1,\dots,D_N$
are bounded pseudoconvex.
\endproclaim

\demo{Proof of Step 2} Let $D_1,\dots,D_N$ be arbitrary pseudoconvex
domains. Let $D_{j,k}\nearrow D_j$, $D_{j,k}\Subset D_j$, where $D_{j,k}$
are pseudoconvex domains with $A_{j,k}:=A_j\cap D_{j,k}\neq\varnothing$.
Observe that all the $A_{j,k}$'s are locally pluriregular. Put
$$
X_k:=\MX(A_{1,k},\dots,A_{N,k};D_{1,k},\dots,D_{N,k})\subset X;
$$
note that $\wdht X_k\nearrow\wdht X$.

Let $f\in\CO_s(X\setminus M)$ be given. By the reduction assumption,
for each $k$ there exists an $\wdht f_k\in\CO(\wdht X_k\setminus M)$ with
$\wdht f_k=f$ on $X_k\setminus M$. By Lemma 6(a), $\wdht f_{k+1}=\wdht f_k$
on $\wdht X_k\setminus M$. Therefore, gluing the $\wdht f_k$'s, we obtain
an $\wdht f\in\CO(\wdht X\setminus M)$ with $\wdht f=f$ on $X\setminus M$.
\qed\enddemo

{\it From now on we assume that $M$ is empty or pure one--codimensional and
$D_1,\dots,D_N$ are bounded pseudoconvex.}

\proclaim{Step 3} To prove the Main Theorem it suffices to consider only the
case $N=2$.
\endproclaim

\proclaim{Remark}
In virtue of Theorem 1, Step 3 finishes the proof of the Main Theorem for
$M=\varnothing$.
\endproclaim

\demo{Proof of Step 3} We proceed by induction on $N\geq2$. Suppose that
the Main Theorem is true for $N-1\geq2$. We have to discuss the case
of an $N$--fold cross $X=\MX(A_1,\dots,A_N;D_1,\dots,D_N)$, where $D_1,\dots,
D_N$ are bounded pseudoconvex.
Let $M\subset\wdht X$ be empty or pure one--codimensional.

Let $f\in\CO_s(X\setminus M)$.

Observe that
$$
X=(Y\times A_N)\cup(\wdht A\times D_N),
$$
where
$$
Y:=\MX(A_1,\dots,A_{N-1};D_1,\dots,D_{N-1}),\qquad
\wdht A:=A_1\tdots A_{N-1}.
$$
We also mention that for any $a_N\in A_N$ we have
$$
\{(z_1,\dots,z_{N-1})\in\CC^{k_1}\tdots\CC^{k_{N-1}}\:
(z_1,\dots,z_{N-1},a_N)\in\wdht X\}=\wdht Y.
$$
Now fix an $a_N\in A_N$ such that
$$
M_{a_N}:=\{(z_1,\dots,z_{N-1})\in\wdht Y\: (z_1,\dots,z_{N-1},a_N)\in M\}
\varsubsetneq\wdht Y;
$$
in particular, $M_{a_N}$ is empty or one--codimensional (in $\wdht Y$).
Recall that $A_1$, \dots, $A_{N-1}$ are locally pluriregular.
By inductive assumption there exists an
$\wdht f_{a_N}\in\CO(\wdht Y\setminus M_{a_N})$ with $\wdht f_{a_N}=
f(\cdot,a_N)$ on $Y\setminus M_{a_N}$.

To continue define the following $2$--fold cross
$$
Z:=\MX(\wdht A,A_N;\wdht Y,D_N).
$$
Notice that $Z$ satisfies all the properties for the case $N=2$:
$\wdht Y, D_N$ are bounded pseudoconvex domains, $\wdht A\subset\wdht Y$,
$A_N\subset D_N$ are locally pluriregular.

In virtue of Lemma 4, we have
$$
\wdht Z=\{(\wdht z,z_N)\in\wdht Y\times D_N\: h_{\wdht A,\wdht Y}
^\ast(\wdht z)+h_{A_N,D_N}^\ast(z_N)<1\}=\wdht X.
$$
Define $\wdtl f\:Z\setminus M\too\CC$,
$$
\wdtl f(z)=\wdtl f(\wdht z,z_N):=\cases \wdht f_{z_N}(\wdht z) &\text{ if }
z\in\wdht Y\times A_N\\
f(z) &\text{ if } z\in\wdht A\times D_N
\endcases.
$$
Obviously, $\wdtl f$ is well-defined and therefore
$\wdtl f\in\CO_s(Z\setminus M)$.

Using the case $N=2$, we find another function $\wdht f
\in\CO(\wdht Z\setminus M)$ with $\wdht f=\wdtl f$ on $Z\setminus M$.
Recall that $\wdht Z=\wdht X$. Hence $\wdht f=f$ on $X\setminus M$.
\qed\enddemo

What remains is to prove the case $N=2$ and $M\neq\varnothing$.
{\it From now on we simplify our notation and we are looking for the following
configuration:

Let $A\subset D\Subset\CC^p$, $B\subset G\Subset\CC^q$, where $D$, $G$ are
bounded pseudoconvex domains, $A, B$ are locally pluriregular.
Put, as always,
$$
X:=\MX(A,B;D,G),\qquad \wdht X:=\{(z,w)\in D\times G\: h_{A,D}^\ast(z)+
h_{B,G}^\ast(w)<1\}.
$$
Moreover, let $M$ be a pure one--codimensional analytic subset of $\wdht X$.}

We like to show that any $f\in\CO_s(X\setminus M)$ extends holomorphically
to $\wdht X\setminus M$.

\proclaim{Step 4} Let $X$, $M$, and $f$ be as above. Let $(D_j)_{j=1}^\infty$,
$(G_j)_{j=1}^\infty$ be sequences of pseudoconvex domains, $D_j\Subset D$,
$G_j\Subset G$, with $D_j\nearrow D$, $G_j\nearrow G$. Moreover, let
$A'\subset A$, $B'\subset B$ be such that $A\setminus A'$,
$B\setminus B'$ are pluripolar, and $A'\cap D_j\neq\varnothing$,
$B'\cap G_j\neq\varnothing$, $j\in\NN$.
For each $j\in\NN$ we assume that:

\halfskip

\noindent for any $(a,b)\in(A'\cap D_j)\times(B'\cap G_j)$ there exist polydiscs
$\Delta_a(r_{a,j})\subset D_j$, $\Delta_b(s_{b,j})\subset G_j$ with
$(\Delta_a(r_{a,j})\times G_j)\cup(D_j\times\Delta_b(s_{b,j}))\subset\wdht X$,
and functions $f_{a,j}\in\CO(\Delta_a(r_{a,j})\times G_j\setminus M)$,
$f^{b,j}\in\CO(D_j\times\Delta_b(s_{b,j})\setminus M)$ such that

$\bullet$ $f_{a,j}=f$ on $(A'\cap\Delta_a(r_{a,j}))\times G_j\setminus M$,

$\bullet$ $f^{b,j}=f$ on $D_j\times(B'\cap\Delta_b(s_{b,j}))\setminus M$.

Then there exists an $\wdht f\in\CO(\wdht X\setminus M)$ with $\wdht f=f$
on $X\setminus M$.
\endproclaim

\demo{Proof of Step 4} Fix a $j\in\NN$. Put
$$
\gather
\wdtl U_j:=\bigcup_{\Sb a\in A'\cap D_j\\ b\in B'\cap G_j\endSb}
(\Delta_a(r_{a,j})\times G_j)\cup(D_j\times\Delta_b(s_{b,j})),\\
X_j:=((A\cap D_j)\times G_j)\cup(D_j\times(B\cap G_j)).
\endgather
$$
Note that
$$
X'_j:=((A'\cap D_j)\times G_j)\cup(D_j\times(B'\cap G_j))\subset\wdtl U_j.
$$
We like to glue the functions $(f_{a,j})_{a\in A'\cap D_j}$ and
$(f^{b,j})_{b\in B'\cap G_j}$ to obtain a global holomorphic function
$f_j$ on $\wdtl U_j\setminus M$:

Let $a\in A'\cap D_j$, $b\in B'\cap G_j$. Observe that
$$
\align
f_{a,j}&=f \;\text{ on } (A'\cap\Delta_a(r_{a,j}))\times G_j\setminus M,\\
f^{b,j}&=f \;\text{ on } D_j\times(B'\cap\Delta_b(s_{b,j}))\setminus M.
\endalign
$$
Thus $f_{a,j}=f^{b,j}$ on $(A'\cap\Delta_a(r_{a,j}))\times
(B'\cap\Delta_b(s_{b,j}))\setminus M$.
Applying Lemma 6(a), we conclude that
$$
f_{a,j}=f^{b,j}\;\text{ on } (\Delta_{a}(r_{a,j})\times
\Delta_{b}(s_{b,j}))\setminus M.
$$

Now let $a', a''\in A'\cap D_j$ be such that
$\Delta_{a'}(r_{a',j})\cap\Delta_{a''}(r_{a'',j})\neq\varnothing$.
Fix a $b\in B'\cap G_j$. We know that $f_{a',j}=f^{b,j}=f_{a'',j}$ on 
$(\Delta_{a'}(r_{a',j})\cap\Delta_{a''}(r_{a'',j}))\times\Delta_b(r_{b,j})
\setminus M$. Hence, by the identity principle, we conclude that $f_{a',j}=f_{a'',j}$ on 
$(\Delta_{a'}(r_{a',j})\cap\Delta_{a''}(r_{a'',j}))\times G_j\setminus M$.

The same argument works for $b', b''\in B'\cap G_j$.

Consequently, we obtain a function $f_j\in\CO(\wdtl U_j\setminus M)$ with
$f_j=f$ on $X'_j\setminus M$.

Let $U_j$ be the connected component of $\wdtl U_j\cap\wdht X'_j$ with
$X'_j\subset U_j$. Thus we have $f_j\in\CO(U_j\setminus M)$ with
$f_j=f$ on $X'_j\setminus M$.

Recall that $X'_j\subset U_j\subset\wdht X'_j$. We claim that
the envelope of holomorphy of $U_j$ coincides with $\wdht X'_j$.
In fact, let $h\in\CO(U_j)$, then $h|_{X'_j}\in\CO_s(X'_j)$. So,
in virtue of Theorem 1, there exists an $\wdht h\in\CO(\wdht X'_j)$ with
$\wdht h=h$ on $X'_j$. Lemma 6(b${}_2$) implies that $\wdht h=h$ on $U_j$.

Applying the Grauert--Remmert theorem
(cf\. \cite{Jar-Pfl 2000}, Th\. 3.4.7), we find a function $\wdht f_j\in
\CO(\wdht X'_j\setminus M)$ with $\wdht f_j=f_j$ on $U_j\setminus M$.
In particular, $\wdht f_j=f$ on $X'_j\setminus M$.

Since $A\setminus A'$, $B\setminus B'$ are pluripolar, we get
$$
\align
\wdht X'_j&=\{(z,w)\in D_j\times G_j\: h_{A'\cap D_j, D_j}^\ast(z)+
h_{B'\cap G_j, G_j}^\ast(w)<1\}\\
&=\{(z,w)\in D_j\times G_j\: h_{A\cap D_j, D_j}^\ast(z)+
h_{B\cap G_j, G_j}^\ast(w)<1\}=\wdht X_j.
\endalign
$$
So, in fact, $\wdht f_j\in\CO(\wdht X_j\setminus M)$. Using Lemma 6(b${}_1$),
we even see that $\wdht f_j=f$ on $X_j\setminus M$.

Observe that $\bigcup_{j=1}^\infty X_j=X$, $\wdht X_j\subset\wdht X_{j+1}$, and
$\bigcup_{j=1}^\infty\wdht X_j=\wdht X$.
Using again Lemma 6(a), by gluing the $\wdht f_j$'s, we get
a function $\wdht f\in\CO(\wdht X\setminus M)$
with $\wdht f=f$ on $X\setminus M$, which finishes the proof of Step 4.
\qed\enddemo

To apply Step 4 we introduce the following condition (*):

Let $\rho>0, 0<r<R$. Put
$$
\Omega:=\Delta_{a_0}(\rho)\times\Delta_{b_0}(R)\subset\CC^p\times\CC^q,
\qquad
\wdtl\Omega:=\Delta_{a_0}(\rho)\times\Delta_{b_0}(r)\subset\CC^p\times\CC^q.
$$
Let $A\subset\Delta_{a_0}(\rho)\subset\CC^p$ be locally pluriregular,
$a_0\in A$, and let $M$ be a pure one--codimensional analytic subset of
$\Omega$ with $M\cap\wdtl\Omega=\varnothing$.
Put $M_a:=\{w\in\Delta_{b_0}(R)\: (a,w)\in M\}$, $a\in A$.

\halfskip

Then we say that the condition (*) holds if:

\hang{\it For any $R'\in(r,R)$ there exists $\rho'\in(0,\rho)$ such that for any
function $f\in\CO(\wdtl\Omega)$ with $f(a,\cdot)\in\CO(\Delta_{b_0}
(R)\setminus M_a)$, $a\in A$,
there exists an extension $\wdht f\in\CO(\Delta_{a_0}(\rho')\times
\Delta_{b_0}(R')\setminus M)$ with $\wdht f=f$ on
$\Delta_{a_0}(\rho')\times\Delta_{b_0}(r)$.}

\proclaim{Step 5} If condition (*) is satisfied, then
the assumptions of Step 4 are fulfilled.
\endproclaim

\demo{Proof of Step 5} Take $X$, $M$, $f\in\CO_s(X\setminus M)$ as
is in Step 4. Define
$$
A':=\{a\in A\: M_a\neq G\},\qquad B':=\{a\in B\: M^b\neq D\},
$$
where $M_a:=\{w\in G\: (a,w)\in M\}$, $M^b:=\{z\in D\: (z,b)\in M\}$.
It clear that $A\setminus A'$, $B\setminus B'$ are pluripolar.

Let $(D_j)_{j=1}^\infty$, $(G_j)_{j=1}^\infty$ be approximation sequences:
$D_j\Subset D_{j+1}\Subset D$, $G_j\Subset G_{j+1}\Subset G$, $D_j\nearrow D$, $G_j\nearrow G$,
$A'\cap D_j\neq\varnothing$, and $B'\cap G_j\neq\varnothing$, $j\in\NN$.

Fix a $j\in\NN$, $a\in A'\cap D_j$ and let $\Omega_j$ be the set of
all $b\in G_{j+1}$ such that there exist a polydisc
$\Delta_{(a,b)}(r_b)\subset D_j\times G_{j+1}$ and a function
$\wdtl f_b\in\CO(\Delta_{(a,b)}(r_b)\setminus M)$ with $\wdtl f_b=f$
on $(A\cap\Delta_a(r_b))\times\Delta_b(r_b)\setminus M$.

It is clear that $\Omega_j$ is open. Observe that $\Omega_j\neq\varnothing$.
Indeed, since $B\cap G_j\setminus M_a\neq\varnothing$,
we find a point $b\in B\cap G_j\setminus M_a$. Therefore there is a polydisc
$\Delta_{(a,b)}(r)\subset D_j\times G_j\setminus M$. Put
$$
Y:=\MX(A\cap\Delta_a(r),B\cap\Delta_b(r);\Delta_a(r),\Delta_b(r)).
$$
By Theorem 1, we find
$r_b\in(0,r)$ and $\wdtl f_b\in\CO(\Delta_{(a,b)}(r_b))$ with
$\wdtl f_b=f$ on $\Delta_{(a,b)}(r_b)\cap Y\supset
(A\cap\Delta_a(r_b))\times\Delta_b(r_b)$.
Consequently, $b\in\Omega_j$.

Moreover, $\Omega_j$ is relatively closed in $G_{j+1}$.
Indeed, let $c$ be an accumulation point of $\Omega_j$ in $G_{j+1}$ and
let $\Delta_c(3R)\subset G_{j+1}$. Take a point $b\in\Omega_j\cap\Delta_c(R)
\setminus M_a$ and let $r\in(0,r_b]$, $r<2R$,
be such that $\Delta_{(a,b)}(r)\cap M=\varnothing$.
Observe that $\wdtl f_b\in\CO(\Delta_{(a,b)}(r))$ and
$\wdtl f_b(z,\cdot)=f(z,\cdot)\in\CO(\Delta_b(2R)\setminus M_z)$ for any
$z\in A\cap\Delta_a(r)$. Hence, by (*) (with $R':=R$), there exists
an extension $\wdht{\wdtl f_b}\in\CO(\Delta_a(\rho')\times\Delta_b(R)\setminus
M)$ ($\rho'\in(0,r)$) such that $\wdht{\wdtl f_b}=\wdtl f_b$ on
$\Delta_{(a,b)}(r)$. Take an $r_c>0$ so small that
$\Delta_{(a,c)}(r_c)\subset\Delta_a(\rho')\times\Delta_b(R)$ and put
$\wdtl f_c:=\wdht{\wdtl f_b}$ on $\Delta_{(a,c)}(r_c)\setminus M$.
Obviously $\wdtl f_c=\wdht{\wdtl f_b}=f$ on $(A\cap\Delta_a(r_c))
\times\Delta_c(r_c)\setminus M$. Hence $c\in\Omega_j$.

Thus $\Omega_j=G_{j+1}$. There exists a finite set $T\subset\overline G_j$
such that
$$
\overline G_j\subset\bigcup_{b\in T}\Delta_b(r_b).
$$
Define $r_{a,j}:=\min\{r_b\: b\in T\}$. Take $b', b''\in T$ with
$\Delta_{b'}(r_{b'})\cap\Delta_{b''}(r_{b''})\neq\varnothing$.
Then $\wdtl f_{b'}=f=\wdtl f_{b''}$ on $(A'\cap\Delta_a(r_{a,j}))
\times(\Delta_{b'}(r_{b'})\cap\Delta_{b''}(r_{b''}))\setminus M$.
Consequently, by Lemma 6(a), $\wdtl f_{b'}=\wdtl f_{b''}$ on
$\Delta_a(r_{a,j})\times(\Delta_{b'}(r_{b'})\cap\Delta_{b''}(r_{b''}))\setminus M$.
In particular, by gluing the functions $(\wdtl f_b)_{b\in T}$, we get a function
$f_{a,j}\in\CO(\Delta_a(r_{a,j})\times G_j\setminus M)$ such that
$f_{a,j}=f$ on $(A'\cap\Delta_a(r_{a,j}))\times G_j\setminus M$.

Changing the role of $z$ and $w$, we get $f^{b,j}$, $b\in B'\cap G_j$.

Thus the assumptions of Step 4 are fulfilled.
\qed\enddemo

It remains to check (*).

\proclaim{Step 6} The condition (*) is always satisfied,
i.e\. the Main Theorem is true.
\endproclaim

\demo{Proof of Step 6} Fix a function $f\in\CO(\wdtl\Omega)$
such that $f(a,\cdot)\in \CO(\Delta_{b_0}(R)\setminus M_a)$ for any
$a\in A$ with $M_a\neq\Delta_{b_0}(R)$. Define
$$
\multline
R_0^\ast:=\sup\{R'\in[r,R)\: \exists_{\rho'\in(0,\rho]}\;\exists_{
\wdht f\in\CO(\Delta_{a_0}(\rho')\times\Delta_{b_0}(R')\setminus M)}\:\\
\wdht f= f\; \text{ on } \Delta_{a_0}(\rho')\times\Delta_{b_0}(r)\}.
\qquad
\endmultline \tag3
$$

It suffices to show that $R_0^\ast=R$.

Suppose that $R_0^\ast<R$. Fix $R_0^\ast<R'_0<R_0<R$ and
choose $R', \rho', \wdht f$ as in (3)
with $R'\in[r,R_0^\ast)$, $\root{q}\of{{R'}^{q-1}R'_0}>R_0^\ast$.

Write $w=(w',w_q)\in\CC^q=\CC^{q-1}\times\CC$. Put
$\wdtl A:=A\cap\Delta_{a_0}(\rho')$.

\halfskip

Let $A'$ denote the set of all $(a,b')\in\wdtl A\times\Delta_{b'_0}(R')$
which satisfy the following condition:

\halfskip

$({}^*_*)$\quad
There exist $R''\in(R_0,R)$, $\delta>0$, $m\in\NN$,
$c_1,\dots,c_m\in\Delta_{b_{0,q}}(R'')$, $\eps>0$,
and holomorphic functions
$\phi_\mu\:\Delta_{(a,b')}(\delta)\too\Delta_{c_\mu}(\eps)$, $\mu=1,\dots,m$,
such that:

$\bullet$ $\Delta_{(a,b')}(\delta)\subset\Delta_{a_0}(\rho')\times
\Delta_{b'_0}(R')$,

$\bullet$ $\Delta_{c_\mu}(\eps)\Subset\Delta_{b_{0,q}}(R'')$, $\mu=1,\dots,m$,

$\bullet$ $\overline\Delta_{c_\mu}(\eps)\cap\overline\Delta_{c_\nu}(\eps)=\varnothing$
for $\mu\neq\nu$, $\mu,\nu=1,\dots,m$,

$\bullet$ $\wdtl H:=\Delta_{b_{0,q}}(R')\cap H\neq\varnothing$,
where $H:=\Delta_{b_{0,q}}(R'')\setminus
\bigcup_{\mu=1}^m\overline\Delta_{c_\mu}(\eps)$,

$\bullet$ $(\Delta_{(a,b')}(\delta)\times\Delta_{b_{0,q}}(R''))\cap M=
\bigcup_{\mu=1}^m\{(z,w',\phi_\mu(z,w'))\: (z,w')\in\Delta_{(a,b')}(\delta)\}$.

\halfskip

For any $(a,b')\in A'$ define a new cross
$$
Y:=\MX((A\cap\Delta_a(\delta))\times\Delta_{b'}(\delta),\wdtl H;
\Delta_{(a,b')}(\delta),H).
$$
Notice that $Y$ does not intersect $M$. In particular, $\wdht f|_Y\in\CO_s(Y)$.
Hence, by Theorem 1, there exists an $\wdht f_1\in\CO(\wdht Y)$ with $\wdht f_1=
\wdht f$ on $Y$. Take $R'''\in(R_0,R'')$, and $\eps''>\eps'>\eps$
($\eps''\approx\eps$), such that

$\bullet$ $\Delta_{c_\mu}(\eps'')\Subset\Delta_{b_{0,q}}(R''')$,
$\mu=1,\dots,m$,

$\bullet$ $\overline\Delta_{c_\mu}(\eps'')\cap\overline\Delta_{c_\nu}(\eps'')=
\varnothing$ for $\mu\neq\nu$, $\mu,\nu=1,\dots,m$.

Then there exists $\delta'\in(0,\delta]$ such that

$\bullet$ $\Delta_{(a,b')}(\delta')\times H'\subset\wdht Y$, where
$H':=\Delta_{b_{0,q}}(R''')\setminus
\bigcup_{\mu=1}^m\overline\Delta_{c_\mu}(\eps')$.

In particular, $\wdht f_1\in\CO(\Delta_{(a,b')}(\delta')\times H')$.

Fix a $\mu\in\{1,\dots,m\}$. Then $\wdht f_1\in\CO(\Delta_{(a,b')}(\delta')\times
(\Delta_{c_\mu}(\eps'')\setminus\overline\Delta_{c_\mu}(\eps')))$ and
$\wdht f_1(z,w',\cdot)\in\CO(\Delta_{c_\mu}(\eps'')\setminus
\{\phi_\mu(z,w')\})$ for any
$(z,w')\in (A\cap\Delta_a(\delta'))\times\Delta_{b'}(\delta')$.
Using the biholomorphic mapping
$$
\align
\Phi_\mu&\:\Delta_{(a,b')}(\delta')\times\CC\too\Delta_{(a,b')}(\delta')\times\CC,
\\
\Phi_\mu&(z,w',w_q):=(z,w',w_q-\phi_\mu(z,w')),
\endalign
$$
we see that the function $g:=\wdht f_1\circ\Phi_\mu^{-1}$ is holomorphic
in $\Delta_{(a,b')}(\delta'')\times
(\Delta_0(\eta'')\setminus\overline\Delta_0(\eta'))$ for some
$\delta''\in(0,\delta']$ and $\eps'<\eta'<\eta''<\eps''$. Moreover,
$g(z,w',\cdot)\in\CO(\Delta_0(\eta'')\setminus\{0\})$ for any
$(z,w')\in (A\cap\Delta_a(\delta''))\times\Delta_{b'}(\delta'')$.
Using Theorem 1 for the cross
$$
\MX((A\cap\Delta_a(\delta''))\times\Delta_{b'}(\delta''),
\Delta_0(\eta'')\setminus\overline\Delta_0(\eta');
\Delta_{(a,b')}(\delta''), \Delta_0(\eta'')\setminus\{0\}),
$$
immediately shows that $g$ extends holomorphically to
$\Delta_{(a,b')}(\delta'')\times(\Delta_0(\eta'')\setminus\{0\})$
(because $h^\ast_{\Delta_0(\eta'')\setminus\{0\},
\Delta_0(\eta'')\setminus\overline\Delta_0(\eta')}\equiv0$).

Transforming the above information back via $\Phi_\mu$ for all $\mu$,
we conclude that the function $\wdht f_1$ extends holomorphically to
$\Delta_{(a,b')}(\delta''')\times\Delta_{b_{0,q}}(R''')\setminus M$ for
some $\delta'''\in(0,\delta'']$; in particular, $\wdht f_1$ extends
holomorphically to
$\Delta_{(a,b')}(\delta''')\times\Delta_{b_{0,q}}(R_0)\setminus M$.

\halfskip

Now we prove that $(\wdtl A\times\Delta_{b'_0}(R'))\setminus A'$ is pluripolar. Write
$$
M\cap(\Delta_{a_0}(\rho')\times\Delta_{b'_0}(R')
\times\Delta_{b_{0,q}}(R))=\bigcup_{\nu=1}^\infty\{\zeta\in P_\nu\: g_\nu(\zeta)=0\},
$$
where $P_\nu\Subset\Delta_{a_0}(\rho')\times\Delta_{b'_0}(R')
\times\Delta_{b_{0,q}}(R)$ is a polydisc and $g_j\in\CO(P_j)$ is a defining
function for $M\cap P_j$; cf\. \cite{Chi 1989}, \S\;2.9. Define
$$
S_\nu:=\{\zeta=(\wdtl\zeta,\zeta_{p+q})\in P_\nu\: g_\nu(\zeta)=
\pder{g_\nu}{\zeta_{p+q}}(\zeta)=0\}
$$
and observe that, by the implicit function theorem, any point from
$$
(\wdtl A\times\Delta_{b'_0}(R'))
\setminus\bigcup_{\nu=1}^\infty\proj_{\wdtl\zeta}(S_\nu)
$$
satisfies (${}^*_*$). It is enough to show that each set $\proj_{\wdtl\zeta}(S_\nu)$
is pluripolar. Fix a $\nu$. Let $S$ be an irreducible component
of $S_\nu$. We have to show that  $\proj_{\wdtl\zeta}(S)$ is pluripolar.
If $S$ has codimension $\geq2$, then $\proj_{\wdtl\zeta}(S)$ is contained in a countable
union of proper analytic sets (cf\. \cite{Chi 1989}, \S\;3.8). Consequently,
$\proj_{\wdtl\zeta}(S)$ is pluripolar. Thus we may assume that $S$ is pure one--codimensional.
The same argument as above shows that
$\proj_{\wdtl\zeta}(\Sing(S))$ is pluripolar.
It remains to prove that $\proj_{\wdtl\zeta}(\Reg(S))$
is pluripolar). Since $g_\nu$ is a defining function,
for any $\zeta\in\Reg(S)$ there exists a $k\in\{1,\dots,p+q-1\}$ such
that $\pder{g_\nu}{\zeta_k}(\zeta)\neq0$. Thus
$$
\Reg(S)=\bigcup_{k=1}^{p+q-1}T_k,
$$
where $T_k:=\{\zeta\in\Reg(S)\:
\pder{g_\nu}{\zeta_k}(\zeta)\neq0\}$. We only need to prove that each set
$\proj_{\wdtl\zeta}(T_k)$ is pluripolar, $k=1,\dots,p+q-1$. Fix a $k$.
To simplify notation, assume that $k=1$. Observe that, by the implicit
function theorem, we can write
$$
T_1=\bigcup_{\ell=1}^\infty\{\zeta\in Q_\ell\: \zeta_1=\psi_\ell(
\zeta_2,\dots,\zeta_{p+q})\},
$$
where $Q_\ell\subset P_\nu$ is a polydisc, $Q_\ell=Q'_{\ell}\times Q''_{\ell}
\subset\CC\times\CC^{p+q-1}$, and $\psi_\ell\:Q''_{\ell}\too Q'_{\ell}$
is holomorphic, $\ell\in\NN$. It suffices to prove that the projection
of each set $T_{1,\ell}:=\{\zeta\in Q_\ell\: \zeta_1=\psi_\ell(
\zeta_2,\dots,\zeta_{p+q})\}$ is pluripolar. Fix an $\ell$. Since
$$
g_\nu(\psi_\ell(\zeta_2,\dots,\zeta_{p+q}),\zeta_2,\dots,\zeta_{p+q})=0,
\quad (\zeta_2,\dots,\zeta_{p+q})\in Q''_{\ell},
$$
we conclude that $\pder{\psi_\ell}{\zeta_{p+q}}\equiv0$ and consequently
$\psi_\ell$ is independent of $\zeta_{p+q}$. Thus
$\proj_{\wdtl\zeta}(T_{1,\ell})=
\{\zeta_1=\psi_\ell(\zeta_2,\dots,\zeta_{p+q-1})\}$ and therefore the
projection is pluripolar. The proof that
$(\wdtl A\times\Delta_{b'_0}(R'))\setminus A'$ is pluripolar
is completed.

\halfskip
Using Step 4, we conclude that $\wdht f$ extends holomorphically to
the domain $\wdht Y\setminus M$, where
$$
\align
\wdht Y:&=\{(z,w',w_q)\in\Delta_{a_0}(\rho')\times\Delta_{b'_0}(R')
\times\Delta_{b_{0,q}}(R_0)\: \\
&\qquad\qquad h^\ast_{A',\Delta_{a_0}(\rho')\times\Delta_{b'_0}(R')}(z,w')+
h^\ast_{\Delta_{b_{0,q}}(R'),\Delta_{b_{0,q}}(R_0)}(w_q)<1\}\\
&=\{(z,w',w_q)\in\Delta_{a_0}(\rho')\times\Delta_{b'_0}(R')
\times\Delta_{b_{0,q}}(R_0)\: \\
&\qquad\qquad h^\ast_{\wdtl A,\Delta_{a_0}(\rho')}(z)+
h^\ast_{\Delta_{b_{0,q}}(R'),\Delta_{b_{0,q}}(R_0)}(w_q)<1\}
\endalign
$$
(here we have used the product property of the relative extremal function).
Since $R'_0<R_0$, we find a $\rho_q\in(0,\rho']$
and a function $\wdtl f_q\in\CO(\Delta_{a_0}(\rho_q)\times\Delta_{b'_0}(R')
\times\Delta_{b_{0,q}}(R'_0)\setminus M)$ such that
$$
\wdtl f_q=\wdht f \;\text{ on } \Delta_{a_0}(\rho_q)\times\Delta_{b_0}(R')
\setminus M.
$$

\halfskip

If $q=1$ we get a contradiction (because $R'_0>R_0^\ast$).

Let $q\geq2$. Repeating the above argument for the coordinates
$w_\nu$, $\nu=1,\dots,q-1$, we find a $\rho_0\in(0,\rho']$
and a function $\wdtl f$ holomorphic in
$$
\Delta_{a_0}(\rho_0)\times\Big(\bigcup_{\nu=1}^q
\Delta_{(b_{0,1},\dots,b_{0,\nu-1})}(R')\times
\Delta_{b_{0,\nu}}(R'_0)\times\Delta_{(b_{0,\nu+1},\dots,b_{0,q})}(R')
\Big)\setminus M
$$
such that $\wdtl f=\wdht f$ on $\Delta_{a_0}(\rho_0)\times
\Delta_{b_0}(R')\setminus M$. Let $\CH$ denote the envelope of holomorphy of
the domain
$$
\bigcup_{\nu=1}^q
\Delta_{(b_{0,1},\dots,b_{0,\nu-1})}(R')\times
\Delta_{b_{0,\nu}}(R'_0)\times\Delta_{(b_{0,\nu+1},\dots,b_{0,q})}(R').
$$
Applying the Grauert--Remmert theorem, we can extend
$\wdtl f$ holomorphically to
$\Delta_{a_0}(\rho_0)\times\CH\setminus M$, i.e\. there exists
an $\wdht{\wdtl f}\in\CO(\Delta_{a_0}(\rho_0)\times\CH\setminus M)$ with
$\wdht{\wdtl f}=f$ on $\Delta_{a_0}(\rho_0)\times\Delta_{b_0}(r)$.
Observe that
$\Delta_{b_0}(\root{q}\of{{R'}^{q-1}R'_0})\subset\CH$. Recall that
$\root{q}\of{{R'}^{q-1}R'_0}>R_0^\ast$; contradiction.
\qed\enddemo

\proclaim{Remark} \rm Notice that the proof of Step 6 shows that
the following stronger version of (*) is true:
Let $\rho>0, 0<r<R$, $\Omega$, $\wdtl\Omega$, $A$, and $a$ be as in (*).
Let $M$ be a pure one--codimensional analytic subset of $\Omega$
(we do not assume that $M\cap\wdtl\Omega=\varnothing$). Then:

\hang{\it For any $R'\in(r,R)$ there exists $\rho'\in(0,\rho)$ such that for any
function $f\in\CO(\wdtl\Omega\setminus M)$ with $f(a,\cdot)\in\CO(\Delta_{b_0}
(R)\setminus M_a)$, $a\in A$,
there exists an extension $\wdht f\in\CO(\Delta_{a_0}(\rho')\times
\Delta_{b_0}(R')\setminus M)$ with $\wdht f=f$ on
$\Delta_{a_0}(\rho')\times\Delta_{b_0}(r)\setminus M$.}
\endproclaim

\halfskip

\noindent{\bf 4. Proof of the Main Theorem in the general case.}
First observe  that the function $\wdht f$ is uniquely
determined (cf\. \S\; 3).

We proceed by induction on $N$.

Let $D_{j,k}\nearrow D_j$, $D_{j,k}\Subset D_{j,k+1}\Subset D_j$, where $D_{j,k}$ are
pseudoconvex domains with $A_{j,k}:=A_j\cap D_{j,k}\neq\varnothing$. Put
$$
X_k:=\MX(A_{1,k},\dots,A_{N,k};D_{1,k},\dots,D_{N,k})\subset X.
$$

\halfskip

It suffices to show that for each $k\in\NN$ the following condition
(${}^{**}_*$) holds.

\halfskip

\noindent(${}^{**}_*$) \quad There exists a domain $U_k$,
$X_k\subset U_k\subset U\cap\wdht X_k$, such that for any
$f\in\CO_s(X\setminus M)$ there exists an $\wdtl f_k\in\CO(U_k\setminus M)$
with $\wdtl f_k|_{X_k\setminus M}=f|_{X_k\setminus M}$.

\halfskip

Indeed, fix a $k\in\NN$ and observe that $\wdht X_k$ is the envelope of
holomorphy of $U_k$ (cf\. the proof of Step 4). Hence,
in virtue of the Dloussky theorem (cf\. \cite{Jar-Pfl 2000}, Th\. 3.4.8, see
also \cite{Por 2001}),
there exists an analytic subset $\wdtl M_k$ of $\wdht X_k$,
$\wdtl M_k\cap U_k\subset M$, such that $\wdht X_k\setminus\wdtl M_k$ is the
envelope of holomorphy of $U_k\setminus M$. In particular, for each
$f\in\CO_s(X\setminus M)$ there exists an
$\overset\approx\to{f_k}\in\CO(\wdht X_k\setminus\wdtl M_k)$ with
$\overset\approx\to{f_k}|_{U_k\setminus M}=\wdtl f_k$. Let
$\CF_k:=\{\overset\approx\to{f_k}\: f\in\CO_s(X\setminus M)\}\subset
\CO(\wdht X_k\setminus\wdtl M_k)$. It is known that (cf\. \cite{Jar-Pfl 2000},
Prop\. 3.4.5) there exists a pure one--codimensional analytic subset
$\wdht M_k\subset\wdht X_k$, $\wdht M_k\subset\wdtl M_k$, such that any point
of $\wdht M_k$ is singular with respect to $\CF_k$, i.e\.

$\bullet$ any function $\overset\approx\to{f_k}$ extends to a function
$\wdht f_k\in\CO(\wdht X_k\setminus\wdht M_k)$ and

$\bullet$ for any $a\in\wdht M_k$ and an open neighborhood $V$ of $a$,
$V\subset\wdht X_k$, there exists an $f\in\CO_s(X\setminus M)$ such that
$\wdht f_k|_{V\setminus\wdht M_k}$ cannot be holomorphically extended to the
whole $V$.

\halfskip

In particular, $\wdht M_{k+1}\cap\wdht X_k=\wdht M_k$. Consequently,
$\wdht M:=\bigcup_{k=1}^\infty\wdht M_k$ is a pure one--codimensional
analytic subset of $\wdht X$, $\wdht M\cap\bigcup_{k=1}^\infty U_k\subset M$,
and for each $f\in\CO_s(X\setminus M)$, the function
$\wdht f:=\bigcup_{k=1}^\infty\wdht f_k$ is holomorphic on
$\wdht X\setminus\wdht M$ with $\wdht f|_{X\setminus M}=f$.

\halfskip

It remains to prove (${}^{**}_*$). Fix a $k\in\NN$. For any $a=(a_1,\dots,a_N)\in
A_{1,k}\tdots A_{N,k}$ let $\rho=\rho_k(a)$ be such that
$\Delta_a(\rho)\subset D_{1,k}\tdots D_{N,k}$.

If $N\geq4$, then we additionally define $(N-2)$--fold crosses
$$
\align
Y_{k,\mu,\nu}:=\MX(&A_{1,k},\dots,A_{\mu-1,k},A_{\mu+1,k},\dots,
A_{\nu-1,k},A_{\nu+1,k},\dots,A_{N,k};\\
&D_{1,k},\dots,D_{\mu-1,k},D_{\mu+1,k},\dots,
D_{\nu-1,k},D_{\nu+1,k},\dots,D_{N,k}),\\
&\hskip200pt 1\leq\mu<\nu\leq N,
\endalign
$$
and we assume that $\rho$ is so small that
$$
\Delta_{(a_1,\dots,a_{\mu-1},a_{\mu+1},\dots, a_{\nu-1},a_{\nu+1},\dots,a_N)}
(\rho)\subset\wdht Y_{k,\mu,\nu},\quad 1\leq\mu<\nu\leq N.
$$

Since
$\{(a_1,\dots,a_{j-1})\}\times\overline D_{j,k+1}\times\{(a_{j+1},\dots,a_N)\}
\Subset U$, we may assume that
$$
\Delta_{(a_1,\dots,a_{j-1})}(\rho)\times D_{j,k+1}\times
\Delta_{(a_{j+1},\dots,a_N)}(\rho)\subset U,\quad j=1,\dots,N.\tag4
$$
We define $N$--fold crosses
$$
\multline
Z_{k,a,j}:=
\MX(A_1\cap\Delta_{a_1}(\rho),\dots,
A_{j-1}\cap\Delta_{a_{j-1}}(\rho),A_{j,k+1},
A_{j+1}\cap\Delta_{a_{j+1}}(\rho),\dots,\\
A_N\cap\Delta_{a_N}(\rho);
\Delta_{a_1}(\rho),\dots,\Delta_{a_{j-1}}(\rho),D_{j,k+1},
\Delta_{a_{j+1}}(\rho),\dots,\Delta_{a_N}(\rho)),\\
 j=1,\dots,N.
\endmultline
$$
Note that $\wdht Z_{k,a,j}\subset U$. Since
$\{(a_1,\dots,a_{j-1})\}\times\overline D_{j,k}\times\{(a_{j+1},\dots,a_N)\}
\Subset\wdht Z_{k,a,j}$, there exists an $r=r_k(a)$, $0<r\leq\rho$, so small
that
$$
V_{k,a,j}:=\Delta_{(a_1,\dots,a_{j-1})}(r)\times D_{j,k}\times
\Delta_{(a_{j+1},\dots,a_N)}(r)\subset\wdht Z_{k,a,j},\quad j=1,\dots,N.
$$
Put
$$
V_k:=\bigcup_{\Sb a\in A_{1,k}\tdots A_{N,k}\\ j\in\{1,\dots,N\}\endSb}
V_{k,a,j}.
$$
Note that $X_k\subset V_k$. Let $U_k$ be the connected component of
$V_k\cap\wdht X_k$ that contains $X_k$.

\halfskip

In view of (4), the Main Theorem with $U=\wdht X$ (which is already
proved in \S\;3) implies that for any $f\in\CO_s(X\setminus M)$
there exists an extension $\wdht f_{k,a,j}\in\CO(\wdht Z_{k,a,j}\setminus M)$
of $f|_{Z_{k,a,j}\setminus M}$. It remains to glue the functions
$$
\wdtl f_{k,a,j}:=\wdht f_{k,a,j}|_{V_{k,a,j}\setminus M},\quad
a\in A_{1,k}\tdots A_{N,k},\; j=1,\dots,N;
$$
then the function $\wdtl f_k:=
\Big(\bigcup_{\Sb a\in A_{1,k}\tdots A_{N,k}\\ j\in\{1,\dots,N\}\endSb}
\wdtl f_{k,a,j}\Big)\Big|_{U_k\setminus M}$ gives the required extension of
$f|_{X_k\setminus M}$.

\halfskip

To check that the gluing process is possible,
let $a, b\in A_{1,k}\tdots A_{N,k}$, $i,j\in\{1,\dots,N\}$ be such that
$V_{k,a,i}\cap V_{k,b,j}\neq\varnothing$. We have the following two cases:

\halfskip

(a) $i\neq j$: We may assume that $i=N-1$, $j=N$. Write
$w=(w',w'')\in\CC^{k_1+\dots+k_{N-2}}\times\CC^{k_{N-1}+k_N}$. Observe that
$$
V_{k,a,N-1}\cap V_{k,b,N}
=\Big(\Delta_{a'}(r_k(a))\cap\Delta_{b'}(r_k(b))\Big)\times
\Delta_{b_{N-1}}(r_k(b))\times\Delta_{a_N}(r_k(a)).
$$

For $c=(c',c'')$, let
$$
\align
M_{c'}&:=\{w''\in\CC^{k_{N-1}+k_N}\:(c',w'')\in M\},\\
M^{c''}&:=\{w'\in\CC^{k_1+\dots+k_{N-2}}\:(w',c'')\in M\};
\endalign
$$
$M_{c'}$ and $M^{c''}$ are analytic subsets of
$$
\align
U_{c'}&:=\{w''\in\CC^{k_{N-1}+k_N}\:(c',w'')\in U\},\\
U^{c''}&:=\{w'\in\CC^{k_1+\dots+k_{N-2}}\:(w',c'')\in U\},
\endalign
$$
respectively.

We consider the following three subcases:

\halfskip

$N=2$: Then $V_{k,a,1}\cap V_{k,b,2}
=\Delta_{b_1}(r_k(b))\times\Delta_{a_2}(r_k(a))$. We know that
$\wdtl f_{k,a,1}=\wdtl f_{k,b,2}$ on the non-pluripolar set
$(A_1\cap\Delta_{b_1}(r_k(b)))\times(A_2\cap\Delta_{a_2}(r_k(a)))\setminus M$.
Hence, by the identity principle, $\wdtl f_{k,a,1}=\wdtl f_{k,b,2}$ on
$V_{k,a,1}\cap V_{k,b,2}\setminus M$.

\halfskip

$N=3$: Then $V_{k,a,2}\cap V_{k,b,3}=
(\Delta_{a_1}(r_k(a))\cap\Delta_{b_1}(r_k(b))\times
\Delta_{b_2}(r_k(b))\times\Delta_{a_3}(r_k(a))$.
Let $C''$ denote the set of all points
$c''\in(A_2\cap\Delta_{b_2}(r_k(b)))\times
(A_3\cap\Delta_{a_3}(r_k(a)))$ such that the set
$M^{c''}$ has codimension $\geq1$ (i.e\. for any $w'\in M^{c''}$ the
codimension of $M^{c''}$ at $w'$ is $\geq1$). Note that $C''$ is
non-pluripolar. We have $\wdtl f_{k,a,2}(\cdot,c'')=f(\cdot,c'')
=\wdtl f_{k,b,3}(\cdot,c'')$ on
$\Delta_{a_1}(r_k(a))\cap\Delta_{b_1}(r_k(b))\setminus M^{c''}$.

Now, let $c'\in\Delta_{a_1}(r_k(a))\cap\Delta_{b_1}(r_k(b))$
be such that the set $M_{c'}$ has codimension $\geq1$.
Then $\wdtl f_{k,a,2}(c',\cdot)=\wdtl f_{k,b,3}(c',\cdot)$
on $C''\setminus M_{c'}$. Consequently, by the identity principle,
$\wdtl f_{k,a,2}(c',\cdot)=\wdtl f_{k,b,3}(c',\cdot)$
on $\Delta_{b_2}(r_k(b))\times\Delta_{a_3}(r_k(a))\setminus M_{c'}$.
Finally, $\wdtl f_{k,a,2}=\wdtl f_{k,b,3}$ on
$V_{k,a,2}\cap V_{k,b,3}\setminus M$.

\halfskip

If $N\in\{2,3\}$, then we jump directly to (b) and we conclude that the Main
Theorem is true for $N\in\{2,3\}$.

\halfskip

$N\geq4$: Here is the only place where the induction over $N$ is used. We
assume that the Main Theorem is true for $N-1\geq3$.

Similarly as in the case $N=3$, let $C''$ denote the set of all points $c''
\in(A_{N-1}\cap\Delta_{b_{N-1}}(r_k(b)))\times
(A_N\cap\Delta_{a_N}(r_k(a)))$ such that the set $M^{c''}$ has codimension
$\geq1$; $C''$ is non-pluripolar. The function $f_{c''}:=
f(\cdot,c'')$ is separately holomorphic on
$Y_{k,N-1,N}\setminus M^{c''}$. By the inductive assumption, the
function $f_{c''}$ extends to a function $\wdht f_{c''}\in
\CO(\wdht Y_{k,N-1,N}\setminus\wdht M(c''))$, where $\wdht M(c'')$ is an
analytic subset of $\wdht Y_{k,N-1,N}$ with $\wdht M(c'')\subset M^{c''}$
in an open neighborhood of $Y_{k,N-1,N}$. Recall that
$$
\Delta_{a'}(r_k(a))\cup\Delta_{b'}(r_k(b))\subset\wdht Y_{k,N-1,N}.
$$
Since $\wdtl f_{k,a,N-1}(\cdot,c'')=f_{c''}$ on
$\Delta_{a'}(r_k(a))\cap Y_{k,N-1,N}\setminus M^{c''}$ and
$\wdtl f_{k,b,N}(\cdot,c'')=f_{c''}$ on
$\Delta_{b'}(r_k(b))\cap Y_{k,N-1,N}\setminus M^{c''}$, we conclude that
$\wdtl f_{k,a,N-1}(\cdot,c'')=\wdht f_{c''}=\wdtl f_{k,b,N}(\cdot,c'')$ on
$\Delta_{a'}(r_k(a))\cap\Delta_{b'}(r_k(b))\setminus M^{c''}$.

Let $c'\in
\Delta_{a'}(r_k(a))\cap\Delta_{b'}(r_k(b))$
be such that the set $M_{c'}$ has codimension $\geq1$.
Then $\wdtl f_{k,a,N-1}(c',\cdot)=\wdtl f_{k,b,N}(c',\cdot)$
on $C''\setminus M_{c'}$. Consequently, by the identity principle,
$\wdtl f_{k,a,N-1}(c',\cdot)=\wdtl f_{k,b,N}(c',\cdot)$ on
$\Delta_{b_{N-1}}(r_k(b))\times\Delta_{a_N}(r_k(a))\setminus M_{c'}$ and,
finally, $\wdtl f_{k,a,N-1}=\wdtl f_{k,b,N}$ on
$V_{k,a,N-1}\cap V_{k,b,N}\setminus M$.

\halfskip

(b) $i=j$: We may assume that $i=j=N$. Observe that
$$
V_{k,a,N}\cap V_{k,b,N}
=\Big(\Delta_{(a_1,\dots,a_{N-1})}(r_k(a))\cap
\Delta_{(b_1,\dots,b_{N-1})}(r_k(b))\Big)\times D_{N,k}.
$$
By (a) we know that:
$$
\alignat2
\wdtl f_{k,a,N}&=\wdtl f_{k,a,N-1} &&\text{ on }
V_{k,a,N}\cap V_{k,a,N-1}\setminus M,\\
\wdtl f_{k,a,N-1}&=\wdtl f_{k,b,N} &&\text{ on }
V_{k,a,N-1}\cap V_{k,b,N}\setminus M.
\endalignat
$$
Hence
$$
\multline
\wdtl f_{k,a,N}=\wdtl f_{k,b,N} \text{ on }
V_{k,a,N}\cap V_{k,a,N-1}\cap V_{k,b,N}\setminus M\\
=\Big(\Delta_{(a_1,\dots,a_{N-1})}(r_k(a))\cap
\Delta_{(b_1,\dots,b_{N-1})}(r_k(b))\Big)\times\Delta_{a_N}(r_k(a))
\setminus M,
\endmultline
$$
and finally, by the identity principle,
$$
\wdtl f_{k,a,N}=\wdtl f_{k,b,N} \text{ on }
V_{k,a,N}\cap V_{k,b,N}\setminus M.
$$

\halfskip

The proof of the Main Theorem is completed.

\halfskip

\noindent{\bf 5. Proof of Theorem 2.}
To prove Theorem 2 we need the following version of the Hartogs theorem
(cf\. \cite{Chi-Sad 1988}, Th\. 1). Let $E$ denote the unit disc.

\proclaim{Theorem 7}
Let $A\subset E^{N-1}$ be locally pluriregular, let
$U\subset E^{N-1}\times\CC$ be an open neighborhood of $A\times\CC$, and let
$M\subset U$ be a relatively closed set such that $M\cap E^N=\varnothing$ and
for any $a\in A$ the fiber
$M_a:=\{w\in\CC\:(a,w)\in M\}$ is polar. Put $X:=\MX(A,E;E^{N-1},\CC)$.
Then  there exists a relatively closed pluripolar
set $S\subset E^{N-1}\times\CC$, $S\cap E^N=\varnothing$,
such that $S_a\subset M_a$, $a\in A$, and
any function $f\in\CO_s(X\setminus M)$ extends
holomorphically to $\wdht X\setminus S=E^{N-1}\times\CC\setminus S$.
\endproclaim

\demo{Proof} Let $\CF:=\CO_s(X\setminus M)$. It is known that each function
$f\in\CF$ has the univalent domain of existence
$G_f\subset E^{N-1}\times\CC$\; \footnote"${}^{(\ddag)}$"{\;We like to thank Professor
Evgeni Chirka for
explaining us some details of the proof of Theorem 1 in \cite{Chi-Sad 1988}.}.
Let $G$ denote the connected component of $\intt\bigcap_{f\in\CF}G_f$ that contains
$E^N$ and let $S:=E^{N-1}\times\CC\setminus G$. It remains to show that
$S$ is pluripolar.

Take $a\in A$ and $b\in\CC\setminus M_a$. Since $M_a$ is polar, there exists
a curve $\gamma\:[0,1]\too\CC\setminus M_a$
such that $\gamma(0)=0$, $\gamma(1)=b$. Take an $\eps>0$ so small that
$$
\Delta_a(\eps)\times(\gamma([0,1])+\Delta_b(\eps))\subset U\setminus M.
$$
Put $V_b:=E\cup(\gamma([0,1])+\Delta_b(\eps))$ and consider the cross
$$
Y:=\MX(A\cap\Delta_a(\eps),E;\Delta_a(\eps),V_b).
$$
Then $f\in\CO_s(Y)$
for any $f\in\CF$. Consequently, by Theorem 1, $\wdht Y\subset G_f$ for any
$f\in\CF$ and hence $\wdht Y\subset G$.
In particular, $\{a\}\times(\CC\setminus M_a)\subset G$.

Thus $S_a\subset M_a$ for all $a\in A$.
Consequently, by Lemma 5 from \cite{Chi-Sad 1988}, $S$ is pluripolar.
\qed\enddemo

\demo{Proof of Theorem 2}
Fix a point $c=(a,b)\in A\times B\setminus M$ and let $\rho>0$
be such that $\Delta_c(\rho)\cap M=\varnothing$. Consider the cross
$$
Y:=\MX(A\cap\Delta_a(\rho),B\cap\Delta_b(\rho);\Delta_a(\rho),\Delta_b(\rho))
$$
and let $0<r\leq\rho$ be such that $\Delta_c(r)\subset\wdht Y$.
Put $A_0:=A\cap\Delta_a(r)$, $B_0:=B\cap\Delta_b(r)$.
By Theorem 1, for any $f\in\CO_s(X\setminus M)$ there exists an
$\wdtl f\in\CO(\Delta_c(r))$ with $\wdtl f=f$ on $A_0\times B_0$.

Define
$$
U:=\CC\times\Delta_b(r),\quad V:=\Delta_a(r)\times\CC.
$$

Observe that the sets
$$
P:=\{z\in\CC\: M_z\text{ is not polar}\},\quad
Q:=\{w\in\CC\: M^w\text{ is not polar}\}
$$
are polar. Indeed,
let $v\in\PSH(\CC^2)$, $v\not\equiv-\infty$, be such that
$M\subset v^{-1}(-\infty)$. Define
$$
u(z):=\sup\{v(z,w)\:w\in\overline E\}, \quad z\in\CC.
$$
Then $P\subset u^{-1}(-\infty)$. It remains to observe that
$u\in\PSH(\CC)$ and $u\not\equiv-\infty$.

\halfskip

Now, by Theorem 7, there exist relatively closed pluripolar sets
$S\subset U\setminus\Delta_c(r)$, $T\subset V\setminus\Delta_c(r)$,
such that for any $f\in\CO_s(X\setminus M)$ there exist
$\overset{\approx}\to{f_1}\in\CO(U\setminus S)$,
$\overset{\approx}\to{f_2}\in\CO(V\setminus T)$  with
$\overset{\approx}\to{f_1}=\overset{\approx}\to{f_2}=\wdtl f$ on
$\Delta_c(r)$. Consequently, the function $\overset{\approx}\to{f}:=
\overset{\approx}\to{f_1}\cup\overset{\approx}\to{f_2}$ is well
defined on $U\cup V\setminus (S\cup T)$.
Note that $R:=S\cup T$ is a relatively closed pluripolar subset of $W:=U\cup V$.
Moreover, $\CC^2$ is the envelope of holomorphy of $W$.

Now, by \cite{Chi 1993}, there exists a closed pluripolar set
$\wdht M\subset\CC^2$ with $\wdht M\cap W\subset R$ such that
$\CC^2\setminus\wdht M$ is the envelope of holomorphy of
$W\setminus R$. Consequently, for each $f\in\CO_s(X\setminus M)$ there
exists an $\wdht f\in\CO(\CC^N\setminus\wdht M)$ with $\wdht f=f$ on
$A_0\times B_0$.

Define
$$
P':=\{z\in\CC\: (M\cup\wdht M)_z\text{ is not polar}\},\quad
Q':=\{w\in\CC\: (M\cup\wdht M)^w\text{ is not polar}\},
$$
$A':=A\setminus P'$, $B':=B\setminus Q'$, $X':=\MX(A',B';\CC,\CC)$.
Now we argue as in the proof of Theorem 1.1
from \cite{Sic 2000} and we conclude that
$\wdht f=f$ on $X'\setminus(M\cup\wdht M)$.
\qed\enddemo

\Refs

{\redefine\bf{\rm}
\widestnumber{\key}{XXXXXXXX}
%----------------------------
\ref
\key Ale-Zer 2001
\by O\. Alehyane, A\. Zeriahi
\paper Une nouvelle version du th\'eor\`eme d'extension de Hartogs pour les
applications s\'epar\'ement holomorphes entre espaces analytiques
\jour Ann\. Polon\. Math\.
\vol 76
\yr 2001
\pages 245--278
\endref
%------------------------------------------
\ref
\key Chi 1989
\by E\. M\. Chirka
\book Complex Analytic Sets
\publ Kluwer Acad\. Publishers
\yr 1989
\endref
%------------------------------------------
\ref
\key Chi 1993
\by E\. M\. Chirka
\paper The extension of pluripolar singularity sets
\jour Proc\. Steklov Inst\. Math\.
\vol 200
\yr 1993
\pages 369--373
\endref
%------------------------------------------
\ref
\key Chi-Sad 1988
\by E\. M\. Chirka, A\. Sadullaev
\paper On continuation of functions with polar singularities
\jour Math\. USSR-Sb\.
\vol 60
\yr 1988
\pages 377--384
\endref
%----------------------------
\ref
\key Jar-Pfl 2000
\by  M\. Jarnicki, P\. Pflug
\book  Extension of Holomorphic Functions
\publ  de Gruyter Expositions in Mathematics 34, Walter de Gruyter
\yr 2000
\endref
%----------------------------
\ref
\key Jar-Pfl 2001
\by  M\. Jarnicki, P\. Pflug
\paper Cross theorem
\jour IMUJ Preprint 2001/08 (to appear in Ann\. Polon\. Math\. (2001))
\endref
%----------------------------
\ref
\key Kli 1991
\by M\. Klimek
\book Pluripotential Theory
\publ Oxford University Press
\yr 1991
\endref
%----------------------------
\ref
\key Ngu 1997
\by Nguyen Thanh Van
\paper Separate analyticity and related subjects
\jour Vietnam J\. Math\.
\vol 25
\yr 1997
\pages 81--90
\endref
%----------------------------
\ref
\key Ngu-Sic 1991
\by Nguyen Thanh Van, J\. Siciak
Fonctions plurisousharmoniques extr\'emales et syst\`emes
doublement orthogonaux de fonctions analytiques
\jour Bull\. Sc\. Math\.
\yr 1991
\vol 115
\pages 235--244
\endref
%----------------------------
\ref
\key Ngu-Zer 1991
\by Nguyen Thanh Van, A\. Zeriahi
\paper  Une extension du th\'eor\`eme de
Hartogs sur les fonctions s\'epar\'ement analytiques
\inbook Analyse Complexe Multivariables, R\'ecents D\`evelopements,
A\. Meril (ed.), EditEl, Rende
\yr 1991
\pages 183--194
\endref
%----------------------------
\ref
\key Ngu-Zer 1995
\by Nguyen Thanh Van, A\. Zeriahi
\paper Syst\`emes doublement othogonaux de fontions holomorphes et applications
\jour Banach Center Publ\.
\vol 31
\yr 1995
\pages 281--297
\endref
%----------------------------
\ref
\key \"Okt 1998
\by O\. \"Oktem
\paper Extension of separately analytic functions and applications to
range characterization of exponential Radon transform
\jour Ann\. Polon\. Math\.
\vol 70
\yr 1998
\pages 195--213
\endref
%----------------------------
\ref
\key \"Okt 1999
\by O\. \"Oktem
\paper Extending separately analytic functions in $\CC^{n+m}$ with
singularities
\inbook Extension of separately analytic functions and applications to
mathematical tomography (Thesis)
\publ Dep\. Math\. Stockholm Univ\.
\yr 1999
\endref
%----------------------------
\ref
\key Por 2001
\by  E\. Porten
\paper On the Hartogs--phenomenon and extension of analytic hypersurfaces
in non-separated Riemann domains
\jour Preprint
\yr 2001
\endref
%----------------------------
\ref
\key Shi 1989
\by B\. Shiffman
\paper On separate analyticity and Hartogs theorem
\jour Indiana Univ\. Math\. J\.
\vol 38
\yr 1989
\pages 943--957
\endref
%----------------------------
\ref
\key Sic 1969
\by J\. Siciak
\paper Separately analytic functions and envelopes
of holomorphy of some lower dimensional subsets of $\CC^n$
\jour Ann\. Polon\. Math\.
\vol 22
\yr 1969--1970
\pages 147--171
\endref
%----------------------------
\ref
\key Sic 1981
\by J\. Siciak
\paper Extremal plurisubharmonic functions in $\CC^N$
\jour Ann\. Polon\. Math\.
\vol 39
\yr 1981
\pages 175--211
\endref
%----------------------------
\ref \key Sic 2000
\by J\. Siciak
\paper Holomorphic functions with singularities on algebraic sets
\jour IMUJ Pre\-print 2000/21 (to appear in Univ\. Iag\. Acta Math\. (2001))
\endref
%----------------------------
\ref
\key Zah 1976
\by V\. P\. Zahariuta
\paper Separately analytic functions, generalizations of
Hartogs theorem, and envelopes of holomorphy
\jour Math\. USSR-Sb\.
\vol 30
\yr 1976
\pages 51--67
\endref
%----------------------------
}
\endRefs

\enddocument